\documentclass[10pt,reqno]{article}
\usepackage{amsmath,amsthm,amsfonts,amssymb,amscd}
\usepackage[dvips]{graphicx}
\usepackage{psfrag}

\input{psfig}

\theoremstyle{plain}
\newtheorem{thm}{Theorem}[section]

\newtheorem{prop}[thm]{Proposition}
\newtheorem{clly}[thm]{Corollary}
\newtheorem{lemma}[thm]{Lemma}

\newtheorem{defi}[thm]{Definition}

\newtheorem{maintheorem}{Theorem}

\title{Transitivity and homoclinic classes for
singular-hyperbolic systems}
\author{C. A. Morales, M. J. Pacifico
\thanks{
This work is partially supported by CNPq, FAPERJ and PRONEX-Dyn. Sys./Brazil.}}

\begin{document}
\maketitle
\begin{abstract}
A
{\em singular hyperbolic set} is
a partially hyperbolic set with singularities
(all hyperbolic) and volume expanding central direction \cite{MPP1}.
We study
connected, singular-hyperbolic, attracting
sets with dense closed orbits
{\em and only one singularity}.
These sets are shown to be transitive for most
$C^r$ flows in the Baire's second category sence.
In general
these sets are shown to be either
transitive or the union of two homoclinic classes.
In the later case we prove the existence of 
finitely many homoclinic classes.
Our results generalize for singular-hyperbolic systems
a well known result for hyperbolic systems
in \cite{N}.
\end{abstract}

\section{Introduction}

A {\em singular hyperbolic set} is
a partially hyperbolic set with singularities
(hyperbolic ones) and volume expanding central direction \cite{MPP1}.
The singular-hyperbolic sets were introduced in \cite{MPP1}
to study robust transitive sets
for three-dimensional flows.
There are some works dealing with such sets
\cite{C,M1,M2,MP1,MP2,MPP2,PPV,Y}.
Some of them explore the similarities
between the
hyperbolic and singular-hyperbolic sets.
We further explore such similarities
motivated by the following well known
fact \cite{N}.

\begin{prop}
\label{hyperbolic}
Let $H$ be a non-trivial hyperbolic set
of a $C^r$ flow $X$ on a closed manifold,
$r\geq 1$.
Suppose that the following properties hold.

\begin{enumerate}
\item
$H$ is connected.
\item
$H$ is attracting, i.e. $H=\cap_{t>0}X_t(U)$ for some
compact neighborhood $U$.
\item
The closed orbits contained in
$H$ are dense in $H$.
\end{enumerate}
Then, $H$ is a homoclinic class.
In particular, $H$ is transitive.
\end{prop}

The obviuos question is whether
the conclusions above hold
replacing hyperbolic by singular-hyperbolic.
This question
has negative answer \cite{MP3}.
Nevertheless we shall prove that
singular-hyperbolic sets
which are
connected, attracting,
has dense closed orbits
and only one singularity are transitive for
most
$C^r$ three-dimensional flows
in the Baire's second category sence.
In general
the considered sets are shown to be either
transitive or the union of two homoclinic classes.
In the later case we prove the existence of 
finitely many homoclinic classes.
Our results generalize Proposition
\ref{hyperbolic} for singular-hyperbolic systems.

Let us state our results in a precise way.
Throughout $M$ denotes a closed three-manifold and
${\cal X}^r$ denotes the space of $C^r$ flows
on $M$ endowed with the
usual $C^r$ topology, $r\geq 1$.
Every flow $X=X_t$ in this work is
orientable, namely it is
tangent to a $C^r$ vector field still denoted
by  $X$.

A subset of ${\cal X}^r$ is
{\em residual} if it is a countable intersection
of open-dense subsets of ${\cal X}^r$.
A property holds for generic flows in ${\cal X}^r$
if it holds for every flow in a residual subset of
${\cal X}^r$.
A compact invariant set $A$ of $X$ is:
{\it attracting} if
$A=\cap_{t>0}X_t(U)$
for some compact neighborhood $U$;
{\em non-trivial} if it is not a closed orbit of $X$;
and 
{\em transitive} if
$A=\omega_X(p)$ for some $p\in A$. Here
$\omega_X(p)$ denotes the set of accumulation points
of the positive orbit of $p$
(this is the $\omega$-limit set of $p$).

A {\em homoclinic class} of $X$
is the closure $H_X(z)$ of the transversal intersection points
of the stable and unstable manifolds of
a hyperbolic periodic orbit $z$ of $X$.
It is known by the Birkhoff-Smale Theorem
that $H_X(z)$ is transitive
with dense periodic orbits.

A compact invariant set
$\Lambda$ of $X$ is {\em partially hyperbolic}
if there are an invariant splitting
$T\Lambda=E^s\oplus E^c$
and positive constants $K,\lambda$ such that:

\begin{enumerate}
\item
{\em $E^s$ is contracting}, namely
$$
\mid\mid DX_t/E^s_x\mid \mid\leq
K e^{-\lambda t},
\,\,\,\,\,\,\forall x\in \Lambda,\,\,\forall
t>0.
$$
\item
{\em $E^s$ dominates $E^c$}, namely
$$
\mid \mid DX_t/E^s_x\mid\mid
\cdot
\mid\mid
DX_{-t}/E^c_{X_{t}(x)}
\mid\mid
\leq K e^{-\lambda t},
\,\,\,\,\,\,\forall x\in \Lambda,\,\,\forall
t>0.
$$
\end{enumerate}

The central direction $E^c$ of
$\Lambda$ is said to be {\em volume expanding}
if the additional condition
$$
\mid J(DX_t/E^c_x)\mid
\geq K e^{\lambda t}
$$
holds $\forall x\in \Lambda$, $\forall t>0$
where
$J(\cdot)$ means the jacobian.

\begin{defi}
\label{DEF}
A partially hyperbolic set is
{\em singular-hyperbolic}
if it has singularities
(all hyperbolic)
and volume expanding central direction.
\end{defi}

The most representative example
of a singular-hyperbolic
set is the geometric Lorenz attractor \cite{ABS,GW}.
Our first result is the following.

\begin{maintheorem}
\label{generico}
The following property holds
for generic flows $X\in
{\cal X}^r$.
Let $\Lambda$ be a singular-hyperbolic set
of $X$. Suppose that the properties below hold.

\begin{enumerate}
\item
$\Lambda$ is connected.
\item
$\Lambda$ is attracting.
\item
The closed orbits contained in $\Lambda$ are dense
in $\Lambda$.
\item
$\Lambda$ has a unique singularity.
\end{enumerate}
Then, $\Lambda$ is transitive.
\end{maintheorem}

Standard $C^1$-topology arguments
\cite{MP1,CMP},
Theorem \ref{generico} and Proposition \ref{hyperbolic}
yield the corollary below.

\begin{clly}
\label{penta}
The following property holds
for generic flows $X\in
{\cal X}^1$.
Let $\Lambda$ be a non-trivial, compact, invariant
set of $X$ containing at most one singularity.
If $\Lambda$ is
connected, attracting and has dense closed orbits
then $\Lambda$
is a homoclinic class of $X$.
\end{clly}

The non-generic case
of Theorem \ref{generico} is considered
below.

\begin{maintheorem}
\label{thA'}
Let $\Lambda$ be a singular-hyperbolic set
of a $C^r$ flow on a closed three-manifold,
$r\geq 1$.
Suppose that the properties below hold.

\begin{enumerate}
\item
$\Lambda$ is connected.
\item
$\Lambda$ is attracting.
\item
The closed orbits contained in $\Lambda$ are dense in $\Lambda$.
\item
$\Lambda$ has a unique singularity.
\end{enumerate}
Then, $\Lambda$ is either

\begin{itemize}
\item
transitive or
\item
non-transitive and
the union
of two homoclinic classes.
\end{itemize}
\end{maintheorem}

The two alternatives of the above theorem
can occur \cite{MP3}.
Theorem \ref{thA'}
implies that a singular-hyperbolic set
of a three-dimensional flow
which is attracting, has dense
closed orbits and only one singularity
splits in a finite (possibly non-disjoint)
union of transitive sets.
This provides a sort of spectral decomposition for
singular-hyperbolic attracting sets
of three-dimensional flows.

Let us make some comments
before the idea of the proof.
The local product structure of hyperbolic sets
is the fundamental tool behind the proof of
Proposition \ref{hyperbolic}.
Indeed this property implies that
hyperbolic sets are
{\em tame}, i.e.
contain finitely
many homoclinic classes \cite[p. 26]{V}.
We cannot use local product structure
to prove Theorem \ref{thA'} since
such a property is not valid
for singular-hyperbolic sets.
In particular, we don't know whether
singular-hyperbolic attracting sets
for three-dimensional flows are tame.
On the other hand, we know that the homoclinic classes
in a singular-hyperbolic
attracting set may be non-disjoint \cite{MP3}.

Let us give the idea of the proofs.
Let $\Lambda$ be a singular-hyperbolic set
of a three-dimensional flow
$X$ satisfying the hypotheses (1)-(4) of Theorem \ref{thA'}.
Let $\sigma$ be the unique singularity
of $\Lambda$.
Then the eigenvalues
$\lambda_1,\lambda_2,\lambda_3$ of
$\sigma$ are real and satisfy
$\lambda_2<\lambda_3<0<\lambda_1$ and
$-\lambda_3<\lambda_1$ for some order (\cite{MPP1}).
It follows that $\sigma$ has
a stable manifold $W^s_X(\sigma)$
(tangent to $\lambda_2,\lambda_3$),
a unstable manifold $W^u_X(\sigma)$
(tangent to $\lambda_1$) and a strong stable manifold
$W^{ss}_X(\sigma)$
(tangent to $\lambda_2$).
Clearly $W^s_X(\sigma)$ is two-dimensional,
$W^{ss}_X(\sigma)$ is one-dimensonal
and $W^{ss}_X(\sigma)\subset W^s_X(\sigma)$.
In particular, $W^s_X(\sigma)\setminus
W^{ss}_X(\sigma)$ has
two connected components
$W^{s,+}_X(\sigma),W^{s,-}_X(\sigma)$.
We define
$H^+_X(\sigma,\Lambda)$ as the
closure of the periodic orbits in $\Lambda$ whose unstable manifold
intersects $W^{s,+}_X(\sigma)$ and
similarly we define $H^-_X(\sigma,\Lambda)$.

In Corollary \ref{cara} we prove
$\Lambda=H^+_X(\sigma,\Lambda)\cup H^-_X(\sigma,\Lambda)$.
In
Theorem \ref{th1} we prove that
if $\Lambda$ is not transitive then
$\omega_X(a)$ is a periodic orbit for all
$a\in W^u_X(\sigma)\setminus\{\sigma\}$.
The proof uses the arguments
in \cite[Theorem 5.2]{MP1} and
Propositions \ref{l2}, \ref{co0}, \ref{co3}.
We shall prove
Theorem \ref{generico} at the end of Section 3
by using
Theorem \ref{th1} and the Kupka-Smale
Theorem \cite{dMP}.
In
Theorem \ref{homoclinic class} we prove
that if $\Lambda$ is not transitive then
both $H^+_X(\sigma,\Lambda)$ and 
$H^-_X(\sigma,\Lambda)$ are homoclinic classes.
Theorem \ref{thA'} easily follows from
Corollary \ref{cara}
and Theorem \ref{homoclinic class}.

The proof
of Theorem \ref{homoclinic class} goes as follows.
By Theorem \ref{th1} we can fix
$a\in W^u_X(\sigma)\setminus\{\sigma\}$
such that $\omega_X(a)=O_0$ for some
periodic orbit $O_0$ of $X$.
In Lemma \ref{c.autovalorreais} we prove
that the expanding eigenvalue
of $O_0$ is positive.
Then $O_0$ divides its unstable manifold
in two components denoted by
$W^{u,+}_X(p_0),W^{u,-}_X(p_0)$
for some $p_0\in O_0$ (Definition \ref{u,+-}).
In Proposition
\ref{homoclinic3} we prove that $H^+_X(\sigma,
\Lambda)=Cl(W^{u,+}_X(p_0))$ and similarly for $-$.
Then, in order to prove Theorem \ref{homoclinic class}
it suffices to prove that $Cl(W^{u,+}_X(p_0))$ is
a homoclinic classs and similarly for $-$.
To prove that $Cl(W^{u,+}_X(p_0))$ is
a homoclinic classs we use 
Propositions \ref{preliminar} and \ref{dense II}.
This finishes the proof of Theorem
\ref{homoclinic class}.
Our methods provide a partial solution
for a question above:

\begin{clly}
\label{galinha}
Singular-hyperbolic sets for three-dimensional
flows which are connected, attracting, {\em non-transitive},
have dense closed orbits and only one singularity
are tame.
\end{clly}

Unfortunately this corollary
applies only in the non-generic case
by Theorem \ref{thA'}.
It would be interesting to
prove our results when more than one
singularity is involved.

\section{Preliminary results}

In what follows $X=X_t$ is a $C^r$ flow
on a closed three-manifold $M$.
The closure of $B\subset M$ is denoted by $Cl(B)$.
If $A$ is a compact invariant set
of $X$ we denote $Sing_X(A)$
the set of singularites of $X$ in $A$.
We denote by
$Per_X(A)$ the union of the
periodic orbits of $X$ in $A$.
A compact invariant set $H$ of
$X$ is {\em hyperbolic}
if the tangent bundle over $H$ has an
invariant decomposition
$E^s\oplus E^X\oplus E^u$
such that $E^s$ is contracting, $E^u$ is expanding
and $E^X$ is generated by the direction of $X$
\cite{PT}.
The Stable Manifold Theory \cite{HPS}
asserts the existence of the stable manifold
$W^s_X(p)$ and the unstable manifold $W^u_X(p)$
associated to $p\in H$.
These manifolds are respectively tangent to the
subspaces
$E^s_p\oplus E^X_p$ and $E^X_p\oplus E^u_p$
of $T_pM$.
In particular, $W^s_X(p)$ and $W^u_X(p)$
are well defined if $p$ belongs to a hyperbolic
periodic orbit of $X$.
An interesting case is when $H$ is {\em saddle-type},
i.e. $dim(E^s)=dim(E^u)=1$.
In this case $W^s_X(p)$ and $W^u_X(p)$ are
two-dimensional submanifolds of $M$
(recall $dim(M)=3$).
The maps
$p\in H\to W^s_X(p)$ and $p\in H\to W^u_X(p)$ are
continuous (in compact parts).
On the other hand,
a compact, singular, invariant set
$\Lambda$ of $X$ is {\em
singular-hyperbolic}
if all its singularities are hyperbolic and
the tangent bundle over $\Lambda$ has an
invariant decomposition
$E^s\oplus E^c$
such that $E^s$ is contracting, $E^s$
dominates $E^c$ and
$E^c$ is volume expanding
(i.e.
the jacobian of
$DX_t/E^c$ grows exponentially
as $t\to \infty$).
See \cite{MPP1} for the precise definition.
Again the Stable Manifold Theory
asserts the existence of the strong stable manifold
$W^{ss}_X(p)$
associated to $p\in \Lambda$.
This manifold is tangent to the subspace
$E^s_p$ of $T_pM$.
For all
$p\in \Lambda$ we
define
$W^s_X(p)=\cup_{t\in I\!\! R}W^{ss}_X(X_t(p))$.
If $p$ is regular (i.e. $X(p)\neq 0$)
then $W^s_X(p)$ is
a well defined two-dimensional submanifold
of $M$.
The map $p\in \Lambda\to W^s_X(p)$
is continuous (in compact parts)
at the regular points $p$ of
$\Lambda$.
The following definition is in \cite{MPP2}.

\begin{defi}
\label{d0}
A singularity $\sigma$ of $X$ is
{\em Lorenz-like} if the eigenvalues
$\lambda_1,\lambda_2,\lambda_3$
of $\sigma$ are real and satisfy
the relation
$$
\lambda_2<\lambda_3<0<-\lambda_3<\lambda_1.
$$
\end{defi}

A Lorenz-like singularity is hyperbolic,
and so, $W^s_X(\sigma)$ and $W^u_X(\sigma)$
do exist.
Moreover, the eigenspace
of $\lambda_2$ is tangent
to a one-dimensional invariant manifold $W^{ss}_X(\sigma)$.
This manifold is called the {\em strong stable
manifold} of $\sigma$.
Clearly $W^{ss}_X(\sigma)$ splits
$W^s_X(\sigma)$ in two connected components.
This fact motivates the following definition.

\begin{figure}[htv] 
\centerline{
\psfig{figure=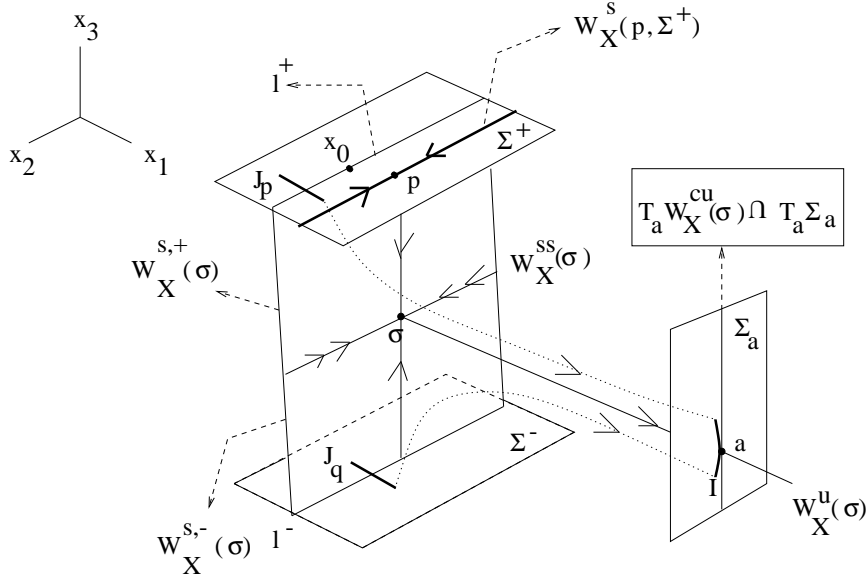,height=3in}}
\caption{\label{f.-2} Singular cross-section.}
\end{figure}

\begin{defi}
\label{W(sing)}
Let $\sigma$ be a Lorenz-like singularity of a
three-dimensional flow $X$.
We denote by
$W^{s,+}_X(\sigma), W^{s,-}_X(\sigma)$
the two connected components
of $W^s_X(\sigma) \setminus
W^{ss}_X(\sigma)$ (see Figure \ref{f.-2}).
\end{defi}

\begin{thm}
\label{th1'}(\cite{MPP2})
Let $\Lambda$ be a singular-hyperbolic set
with dense periodic orbits. Then,
every $\sigma\in Sing_X(\Lambda)$ is Lorenz-like
and satisfies
$\Lambda\cap W^{ss}_X(\sigma)=\{\sigma\}$.
\end{thm}

\begin{lemma}(\cite{MPP2})
\label{lhyp}
A
compact invariant set without singularities
of a singular-hyperbolic set of a three-dimensional flow
is hyperbolic saddle-type.
\end{lemma}

\begin{lemma}
\label{acc->int}
Let $\Lambda$ be a connected, singular-hyperbolic,
attracting set with dense periodic orbits of $X$.
For every $p\in Per_X(\Lambda)$ there is
$\sigma\in Sing_X(\Lambda)$ such that
$$
W^u_X(p)\cap W^s_X(\sigma)\neq \emptyset.
$$
\end{lemma}

\begin{proof}
Fix $p\in Per_X(\Lambda)$.
We claim that $Cl(W^u_X(p))$
is singular.
Indeed, suppose by contradiction that
it is not.
Then,
$Cl(W^u_X(p))$ is hyperbolic saddle
by Lemma \ref{lhyp}.
Using $dim(W^u_X(p))=2$
one can
prove that $Cl(W^u_X(p))$ is attracting, and so,
it is open and closed
in $\Lambda$.
Since $\Lambda$ is connected we conclude that
$\Lambda=Cl(W^u_X(p))$. And this is a
contradiction since
$\Lambda$ contains singularities and
$Cl(W^u_X(p))$ does not.
This contradiction proves the claim.
If $\sigma\in Sing_X(Cl(W^u_X(p)))$ we
can prove
$W^u_X(p)\cap W^s_X(\sigma)\neq\emptyset$
as in the proof of
\cite[Theorem 4.1 p. 365]{MP1}.
\end{proof}

Lemma \ref{acc->int}
and Theorem \ref{th1'}
motivate the following definition.

\begin{defi}
\label{H+}
If $\Lambda$ is a singular-hyperbolic set
of $X$ and $\sigma\in
Sing_X(\Lambda)$ is
Lorenz-like we define

\begin{itemize}
\item
$
P^+_X(\sigma,\Lambda)=\{p\in Per_X(\Lambda):
W^u_X(p)\cap W^{s,+}_X(\sigma)\neq\emptyset\}$.
\item
$
P^-_X(\sigma,\Lambda)=\{p\in Per_X(\Lambda):
W^u_X(p)\cap W^{s,-}_X(\sigma)\neq\emptyset\}$.
\item
$
H^+_X(\sigma,\Lambda)=Cl(P^+_X(\sigma,\Lambda))$.
\item
$
H^-_X(\sigma,\Lambda)=Cl(P^-_X(\sigma,\Lambda))$.
\end{itemize}
\end{defi}

Lemma \ref{acc->int} and Theorem \ref{th1'}
easily imply the following.

\begin{clly}
\label{cara}
Let $\Lambda$ be a connected, singular-hyperbolic,
attracting set with dense periodic orbits
and only one singularity $\sigma$.
Then,
$$
\Lambda=H^+_X(\sigma,\Lambda)\cup
H^-_X(\sigma,\Lambda).
$$
\end{clly}

We finish this section with two technical
(but straightforward) lemmas to be used in the
sequel.
First we state some short definitions.
If  $S$ is a submanifold we shall denote
by $T_xS$ the tangent space at $x\in S$.
By {\em cross-section} of $X$ we mean
a compact submanifold $\Sigma$
transverse to $X$ and diffeomorphic to the
two-dimensional square $[0,1]^2$.
If $\Lambda$ is a singular-hyperbolic set
of $X$ and $x\in \Sigma\cap \Lambda$,
then $x$ is regular and so
$W^s_X(x)$ is a two-dimensional submanifold
transverse to $\Sigma$.
In this case we define
by $W^s_X(x,\Sigma)$ the
connected component of $W^s_X(x)\cap \Sigma$
containing $x$.
We shall be interested in a
special cross-section described as follows.
Let $\Lambda$ be a singular-hyperbolic set
of a three-dimensional flow $X$
and $\sigma\in Sing_X(\Lambda)$.
Suppose that the closed orbits contained in $\Lambda$ are dense in $\Lambda$.
Then $\sigma$ is Lorenz-like
by Theorem \ref{th1'}.
It is then possible to describe the
flow of $X$ close to $\sigma$
by using the Grobman-Hartman Theorem
\cite{dMP}.
Indeed, we can assume that
the flow of $X$ around
$\sigma$ is
the linear flow
given by
$$
X=\lambda_1\partial_{x_1}+\lambda_2\partial_{x_2}+
\lambda_3\partial_{x_3}
$$
in a suitable coordinate system
$(x_1,x_2,x_3)\in [-1,1]^3$ around
$\sigma=(0,0,0)$.
A cross-section $\Sigma$ of $X$ is
{\em singular}
if $\Sigma$ corresponds to the submanifolds
$\Sigma^+=\{z=1\}$ or
$\Sigma^-=\{z=-1\}$ in the coordinate system
$(x_1,x_2,x_3)$.
We denote by $l^+$ and $l^-$
the curves in $\Sigma^+, \Sigma^-$
intersecting to $\{y=0\}$.
Note that $l^+,l^-$ are
contained in $W^{s,+}_X(\sigma),W^{s,-}_X(\sigma)$ respectively.
With these notation
we have the following straighforward lemma.

\vspace{10pt}

\begin{lemma}
\label{sigma}
Let $\Lambda$ a singular-hyperbolic set
with dense periodic orbits of a three-dimensional flow $X$ and $\sigma\in Sing_X(\Lambda)$
be fixed. Then,
there are singular
cross-sections $\Sigma^+,\Sigma^-$
as above such that
every orbit of $\Lambda$ passing close to
some point in $W^{s,+}_X(\sigma)$
(resp. $W^{s,-}_X(\sigma)$)
intersects $\Sigma^+$ (resp. $\Sigma^-$).
If $p\in \Lambda\cap\Sigma^+$ is close to $l^+$,
then
$W^s_X(p,\Sigma^+)$ is a vertical curve
crossing $\Sigma^+$ as indicated in Figure
\ref{f.-2}.
If $p\in Per_X(\Lambda)$ and
$W^u_X(p)\cap W^{s,+}_X(\sigma)\neq\emptyset$, then
$W^u_X(p)$ contains an interval
$J=J_p$ intersecting $l^+$ transversally.
Similarly replacing $+$ by $-$.
\end{lemma}

\begin{proof}
Let $F$ be a fundamental
domain of $W^s_X(\sigma)$
(see \cite{dMP}). Denote $F^+=F\cap W^{s,+}_X(\sigma)$
and $F^-=F\cap W^{s,-}_X(\sigma)$.
By Theorem \ref{th1'} we have
$(\Lambda\cap F)\cap W^{ss}_X(\sigma)=\emptyset$, and so, $\Lambda\cap F=
(\Lambda\cap F^+)\cup (\Lambda\cap F^-)$.
On the other hand, both $\Lambda\cap F^+$ and $\Lambda\cap
F^-$ are compact and far from $W^{ss}_X(\sigma)$.
If follows from the Strong $\lambda$-lemma
\cite{D} that there
are a singular cross-sections
$\Sigma^+,\Sigma^-$ arbitrarily close to $\sigma$
as above
such that
every orbit in
$\Lambda\cap F^+$
(resp. $\Lambda\cap
F^-$) has a point in
$l^+$ (resp. $l^-$). Then,
since $F$ is a fundamental domain,
every orbit in $\Lambda\cap W^{s,+}_X(\sigma)$
intersects $l^+$ (and similarly for
$-$). This implies that every orbit
of $\Lambda$ passing close to some point
in $W^{s,+}_X(\sigma)$ intersects
$\Sigma^+$ (and similarly for $-$). 
Now let $p\in \Lambda\cap \Sigma^+$ be close to
$l^+$.
Since the strong unstable manifolds
$W^{ss}_X(x)$ have size bounded away from zero
($x\in \Lambda$) we have that if
$\Sigma^+$ is close to $\sigma$, then
$W^s_X(x_0,\Sigma^+)=l^+$
($\forall x_0\in \Lambda\cap l^+$).
The continuity of the set-valued map
$p\in \Lambda\to W^s_X(p)$ at the regular points
implies that if
$p\in \Lambda\cap \Sigma^+$ is close
to some point
$x_0\in l^+$ then
$W^s_X(p,\Sigma^+)$ is close to $W^s_X(x_0,\Sigma^+)$.
It follows that $W^s_X(p,\Sigma^+)$
is a vertical curve as desired.
Finally let $p\in Per_X(\Lambda)$ be fixed.
We have that
$W^u_X(p)\subset \Lambda$
since $\Lambda$ is attracting, and so,
$W^u_X(p)\cap W^{s,+}_X(\sigma)\subset
\Lambda\cap W^{s,+}_X(\sigma)$.
As already mentioned every orbit
in $\Lambda\cap W^{s,+}_X(\sigma)$
intersects $l^+$.
In particular,
$W^u_X(p)\cap W^{s,+}_X(\sigma)\neq\emptyset$
implies $W^u_X(p)\cap l^+\neq\emptyset$.
Then, we can construct $J_p$ by
flow-projecting to $\Sigma^+$
nearby $W^u_X(p)\cap l^+$.
\end{proof}

\begin{lemma}
\label{lema*}
Let $\Lambda$ a singular-hyperbolic set
with dense periodic orbits of a three-dimensional
flow $X$ and $\sigma\in Sing_X(\Lambda)$ be fixed.
Let $p,q\in Per_X(\Lambda)$ be such that
$W^u_X(p)\cap W^{s,+}_X(\sigma)\neq\emptyset$
and
$W^u_X(q)\cap W^{s,-}_X(\sigma)\neq\emptyset$.
Let $a\in W^u_X(\sigma)\setminus \{\sigma\}$ and
$\Sigma_a$ be a cross-section of $X$ containing $a$.
Then, there is an open interval $I=I_a\subset \Sigma_a$
containing $a$ such that the properties below hold.

\begin{enumerate}
\item
$I\setminus \{a\}$ is formed by
two intervals, one of them contained in $W^u_X(p)$ and the
other one contained in $W^u_X(q)$.
In particular $I\subset \Lambda$.
\item
$T_xI=E^c_x\cap T_x\Sigma_a$
for all $x\in I$.
\end{enumerate}
\end{lemma}

\begin{proof}
Let $E^s\oplus E^c$ be the singular-hyperbolic splitting
of $\Lambda$.
As the periodic orbits in $\Lambda$ are dense in $\Lambda$ we have that
$X(x)\subset E^c_x$ and $dim(E^c_x)=2$ for all
$x\in \Lambda$.
By Theorem \ref{th1'} we have that $\sigma$ is Lorenz-like (recall Definition \ref{d0}).
By the Invariant Manifold Theory
\cite{HPS}
there is a centre-unstable manifold
$W^{cu}_X(\sigma)$ at $\sigma$.
This manifold is flow-invariant
and tangent to the eigenspace associated to
the eigenvalues
$\lambda_1,\lambda_3$ of $\sigma$.
The centre-unstable manifold is not unique
but all of them satisfy
$W^u_X(\sigma)\subset W^{cu}_X(\sigma)$.
As $E^s$ dominates
$E^c$ we can easily prove that
$T_aW^{cu}_X(\sigma)=E^c_a$.
In particular, one has
$$
T_aW^{cu}_X(\sigma)\cap T_a\Sigma_a=E^c_a\cap T_a\Sigma_a.
$$
See Figure \ref{f.-2}.
Let $\Sigma^+,\Sigma^-$ be the cross-sections given
by Lemma \ref{sigma}.
As
$W^u_X(p)\cap W^{s,+}_X(\sigma)\neq\emptyset$
and
$W^u_X(q)\cap W^{s,-}_X(\sigma)\neq\emptyset$
we can also fix the intervals
$J_p,J_q$ in $W^u_X(p)\cap \Sigma^+,
W^u_X(q)\cap \Sigma^-$ given in Lemma \ref{sigma}.
As $\Lambda$ is attracting we have that
$J_p$ and $J_q$ are both contained in $\Lambda$.
The dominance condition of the splitting
implies $T_xJ_p\subset E^c_x\cap T_x\Sigma^+$ and
$T_{x'}J_q\subset E^c_{x'}\cap \Sigma^-$
for all $x,x'$ in $J_p,J_q$ respectively.
It follows from Lemma \ref{sigma} that
none of the curves $J_p,J_q$ intersects
$W^{ss}_X(\sigma)$.
In addition,
$J_p$ intersects transversally the curve
$l=l^+$ in Lemma \ref{sigma} and similarly for $J_q$.
It follows from the Strong $\lambda$-lemma
\cite{D} that the flow of $X$ carries
some interval
in $J_p$ into an open interval
$I^+\subset \Sigma_a$,
and some interval in $J_q$ into an open interval
$I^-\subset \Sigma_a$,
such that the union $I=I^+\cup I^-\cup\{a\}$ satisfies:
$I$ is an open interval containing $a$ and
$$
T_aI=T_a W^{cu}_X(\sigma)\cap T_a\Sigma_a.
$$
Then,
$$
T_aI=E^c_a\cap T_a\Sigma_a.
$$
As $I\setminus \{a\}=I^+\cup I^-$
we have the property (1) of the lemma
by the invariance of the unstable manifolds.
That $x\in I\setminus \{a\}$ satisfies
the property (2) of the lemma is clear.
That $x=a$ also satisfies the property (2) follows from
$T_aI=E^c_a\cap T_a\Sigma_a$. The lemma is proved.
\end{proof}



\section{Non-transitive singular-hyperbolic sets}

In this section we study
non-transitive singular-hyperbolic sets
for $C^r$ flows on closed three-manifolds
$r\geq 1$.
The result is the following.

\begin{thm}
\label{th1}
Let $\Lambda$ be a singular-hyperbolic set
of a $C^r$ flow $X$ on a closed three-manifold.
Suppose that the properties below hold.

\begin{enumerate}
\item
$\Lambda$ is connected.
\item
$\Lambda$ is attracting.
\item
The closed orbits contained in $\Lambda$ are dense
in $\Lambda$.
\item
$\Lambda$ has a unique singularity $\sigma$.
\item
$\Lambda$ is not transitive.
\end{enumerate}
Then,
$\omega_X(a)$ is a periodic orbit of $X$
for every $a\in W^u_X(\sigma)\setminus \{\sigma\}$.
\end{thm}

Before the proof
we state some preliminars.
First of all we observe that the hypothesis (3)
of the theorem implies

\begin{description}
\item{(H1)}
$\Lambda=Cl(Per_X(\Lambda))$.
\end{description}

To prove transitivity we shall use the
following criterium due to Birkhoff.

\begin{lemma}
\label{trans}
Let $T$ be a compact, invariant set of $X$
such that for all open sets
$U,V$ intersecting $T$ there is
$s>0$ such that
$X_s(U\cap T)\cap V\neq\emptyset$.
Then $T$ is transitive.
\end{lemma}

Birkhoff' criterium will be applied together
with the following lemma.

\begin{lemma}
\label{l1}
Let $\Lambda$ be a singular-hyperbolic set
of $X$ satisfying (1)-(4)
of Theorem \ref{th1}.
Let $U,V$ be open sets,
$p\in U\cap Per_X(\Lambda)$ and $q\in
V\cap Per_X(\Lambda)$.
If $W^u_X(p)\cap W^{s,+}_X(\sigma)\neq\emptyset$
and $W^u_X(q)\cap W^{s,+}_X(\sigma)\neq\emptyset$,
then there are $z\in W^u_X(p)$ arbitrarily close to $W^u_X(p)\cap W^{s,+}_X(\sigma)$
and $t>0$ such that
$X_t(z)\in V$.
A similar result holds replacing $+$ by $-$.
\end{lemma}

\begin{figure}[htv] 
\centerline{
\psfig{figure=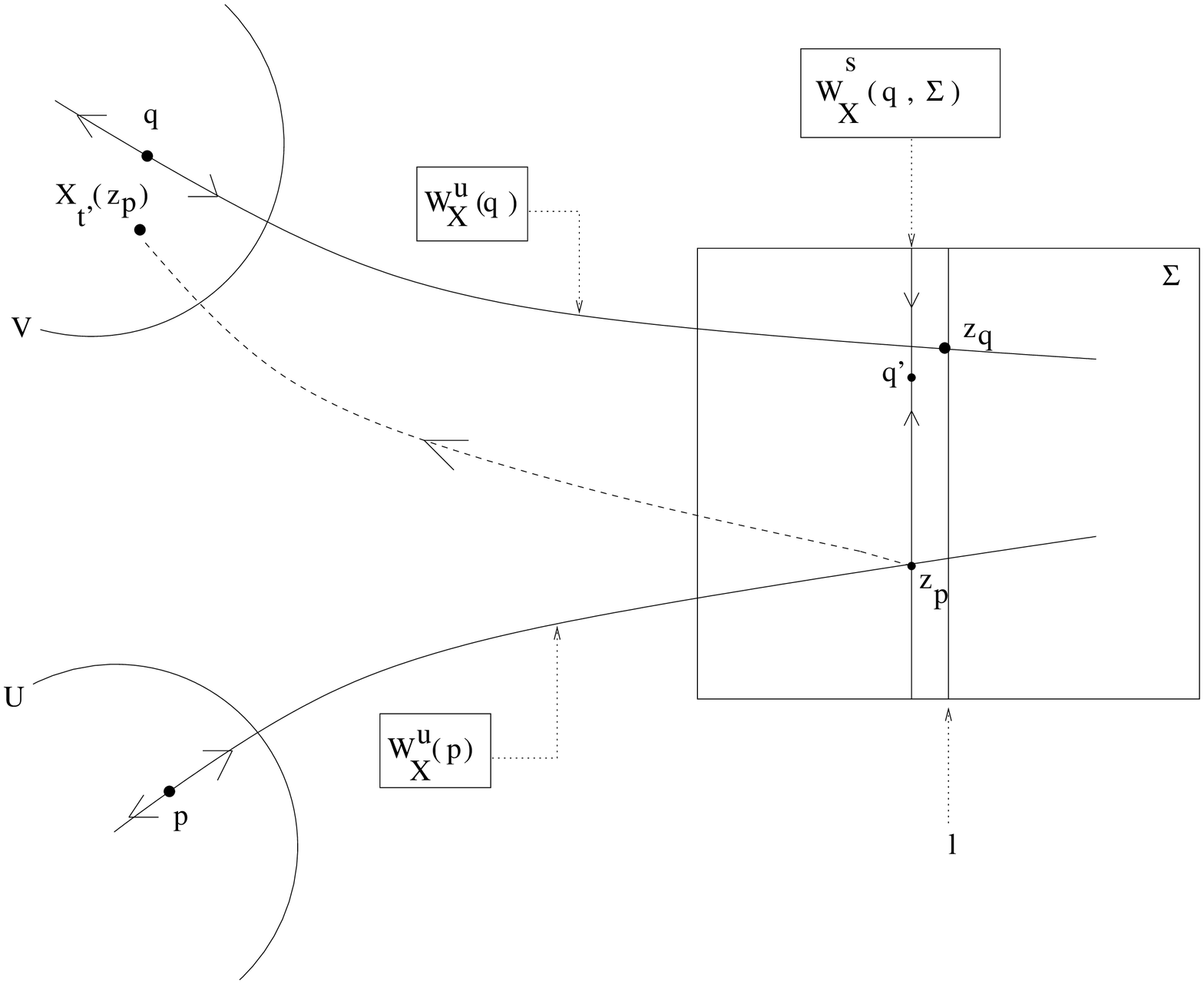,height=4in}}
\caption{\label{f.3} Proof of Lemma \ref{l1}.}
\end{figure}

\begin{proof}
Let $p,q, U, V$ as in the statement.
Fix a cross section $\Sigma=\Sigma^+$ through $W^{s,+}_X(\sigma)$
so that every orbit of $\Lambda$ passing close to
some point in $W^{s,+}_X(\sigma)$ intersects $\Sigma$
(see Lemma \ref{sigma}).
By assumption
$W^u_X(p)\cap W^{s,+}_X(\sigma)\neq\emptyset$
and $W^u_X(q)\cap W^{s,+}_X(\sigma)\neq\emptyset$.
Note that these intersections are transversal.
Because
$W^u_X(q)\cap W^{s,+}_X(\sigma)\neq\emptyset$
there is
$z_q\in W^u_X(q)\cap (W^{s,+}_X(\sigma)\cap \Sigma)$.
By (H1) there is
$q'\in Per_X(\Lambda)\cap \Sigma$ nearby
$z_q$. By Lemma \ref{sigma} the
manifolds $W^s_X(x,\Sigma)$ are well defined vertical
curves,
$\forall x\in \Lambda\cap \Sigma$.
Choosing $q'$ nearby $z_q$,
we have that $W^s_X(q',\Sigma)$ crosses
$\Sigma$ as in Figure \ref{f.3}.
In particular,
as $W^u_X(p)\cap W^{s,+}_X(\sigma)\neq\emptyset$,
there is $z_p\in
W^s_X(q',\Sigma)\cap W^u_X(p)$.
Since $q'$ is close to
$z_q\in W^u_X(q)$ the negative
orbit of $q'$ meets $V$.
Then, as $q'$ is periodic,
we have that the
positive orbit of
$q'$ meets $V$ as well.
As $z_p\in W^s_X(q')$ and $q'\in Per_X(\Lambda)$
there is $t'>0$ such that
$X_{t'}(z_p)\in V$.
The result follows
choosing $z=z_p$ and $t=t'$. 
\end{proof}

The next three propositions give
sufficient conditions for
the transitivity of singular-hyperbolic
sets satisfying (1)-(4) of
Theorem \ref{th1}.

\begin{prop}
\label{l2}
Let $\Lambda$ be a singular-hyperbolic set
of $X$ satisfying (1)-(4)
of Theorem \ref{th1}.
Suppose that
there is a sequence $p_n\in Per_X(\Lambda)$ converging
to some point in $W^{s,+}_X(\sigma)$ such that
$
W^u_X(p_n)\cap W^{s,-}_X(\sigma)\neq\emptyset$,
$\forall n$.
Then $\Lambda$ is transitive.
Similarly
interchanging the roles of $+$ and $-$.
\end{prop}

\begin{proof}
To prove that $\Lambda$
is transitive
we only need to prove by
Birkohff's criterium that
$\forall U,V$ open sets intersecting $\Lambda$
there is $s>0$ such that
$X_s(U\cap \Lambda)\cap V\neq\emptyset$.
For this we proceed as follows:
By (H1) there are
$p\in Per_X(\Lambda)\cap U$
and $q\in Per_X(\Lambda)\cap V$.
Next we proceed according to the following cases:
\item
{\em Case I:
$W^u_X(p)\cap W^{s,+}_X(\sigma)\neq
\emptyset$ and $W^u_X(q)\cap W^{s,+}_X(\sigma)
\neq\emptyset$}.
In this case
Lemma \ref{l1} implies that
there are
$z\in W^u_X(p)$ and $t>0$ such that
$X_t(z)\in V$.
Since $z\in W^u_X(p)$ there is
$t'>0$ such that
$X_{-t'}(z)\in U$.
As $\Lambda$ is attracting
we have $W^u_X(p)\subset \Lambda$,
and so, $z\in \Lambda$.
Choose
$s=t+t'$ and $w=X_{-t'}(z)$.
Then, $s>0$ and
$w\in U\cap \Lambda\cap X_{-s}(V)$.
It follows that $X_s(U\cap \Lambda) \cap V
\neq\emptyset$ for the chosen $s>0$ proving
the result.
\item
{\em Case II:
$W^u_X(p)\cap W^{s,-}_X(\sigma)\neq
\emptyset$ and $W^u_X(q)\cap W^{s,-}_X(\sigma)
\neq\emptyset$}.
In this case the proof is similar
to the one in Case I replacing
$+$ by $-$.
\item
{\em Case III: $W^u_X(p)\cap W^{s,+}_X(\sigma)\neq
\emptyset$ and $W^u_X(q)\cap W^{s,-}_X(\sigma)
\neq\emptyset$}.
Fix a cross section $\Sigma=\Sigma^+$
through $W^{s,+}_X(\sigma)$
and a open arc $J\subset \Sigma\cap
W^u_X(p)$ intersecting $W^{s,+}_X(\sigma)$
transversally (see Lemma \ref{sigma}).
By Lemma \ref{sigma} we can assume that
$p_n\in \Sigma$
for every $n$.
Because
the direction $E^s$ of $\Lambda$
is contracting
the size of $W^s_X(p_n)$ is uniformly
bounded away from zero.
It follows that there is $n$ large
so that $J$
intersects $W^s_X(p_n)$ transversally.
Applying the Inclination Lemma \cite{dMP} to the saturated
of $J\subset W^u_X(p)$,
and the assumption
$
W^u_X(p_n)\cap W^{s,-}_X(\sigma)\neq\emptyset
$,
we conclude that
$W^u_X(p)\cap W^{s,-}_X(\sigma)\neq\emptyset$.
This implies that Case II holds
proving the result.
\item
{\em Case IV: $W^u_X(p)\cap W^{s,-}_X(\sigma)\neq
\emptyset$ and $W^u_X(q)\cap W^{s,+}_X(\sigma)
\neq\emptyset$}.
In this case we interchange the roles
of $p$ and $q$ in the previous
Case III to prove
$W^u_X(q)\cap W^{s,-}_X(\sigma)\neq\emptyset$.
This implies that Case II holds proving the result.
\end{proof}

\begin{prop}
\label{co0}
Let $\Lambda$ be a singular-hyperbolic set
of $X$ satisfying (1)-(4)
of Theorem \ref{th1}.
Suppose that
the following property hold :
If $p,q\in Per_X(\Lambda)$ then
either
\begin{enumerate}
\item
$W^u_X(p)\cap W^{s,+}_X(\sigma)\neq\emptyset$
and $W^u_X(q)\cap W^{s,+}_X(\sigma)\neq\emptyset$ or
\item $W^u_X(p)\cap W^{s,-}_X(\sigma)\neq\emptyset$
and $W^u_X(q)\cap W^{s,-}_X(\sigma)\neq\emptyset$.
\end{enumerate}

Then, $\Lambda$ is transitive.
\end{prop}

\begin{proof}
Again, by Birkhoff's criterium,
we only need to prove that
$\forall U,V$ open sets intersecting $\Lambda$
there is $s>0$ such that
$X_s(U\cap \Lambda)\cap V\neq\emptyset$.
For this we proceed as follows:
By (H1) there are
$p\in Per_X(\Lambda)\cap U$
and $q\in Per_X(\Lambda)\cap V$.
First suppose that
the alternative (1)
holds. Then, by Lemma \ref{l1},
there are $z\in W^u_X(p)$ and $t>0$ such that
$X_t(z)\in V$.
As $z\in W^u_X(p)$ we have that
$w=X_{-t'}(z)\in U$ for some $t'>0$.
As $\Lambda$ is an attracting set
we have that $w\in \Lambda$.
If $s=t+t'>0$ we conclude
that
$w\in (U\cap\Lambda)\cap X_{-s}(V)$ and so
$X_s(U\cap \Lambda)\cap V\neq\emptyset$.
If the alternative (2) of the corollary holds
we can find
$s>0$ such that $X_s(U\cap \Lambda)\cap V\neq\emptyset$
in a similar way (replacing $+$ by $-$).
The proof follows.
\end{proof}

The following corollary
illustrates how
Proposition \ref{co0} will be used.
The corollary will not be used to prove
the results in this paper.
Its statement uses the following
short definition.
Let $\Lambda$ a compact invariant set
of a three-dimensional flow $X$ and
$\sigma\in \Lambda$ a Lorenz-like singularity
of $X$.
We say that $\Lambda$ {\em accumulated $\sigma$
by one side} if either
$\Lambda\cap W^{s,+}_X(\sigma)$
or $\Lambda\cap W^{s,-}_X(\sigma)$ is $\emptyset$
(otherwise $\Lambda$ accumulates
$\sigma$ in both sides).
A typical example of a singular-hyperbolic, attracting set containing a Lorenz-like singularity
accumulated by one side is the
geometric Lorenz attractor \cite{ABS,GW}.
It is not difficult to construct
singular-hyperbolic attractors
with singularities accumulated in both sides.
The corollary is the following.

\begin{clly}
\label{inutil}
Let $\Lambda$ be a singular-hyperbolic set
of $X$ satisfying (1)-(4)
of Theorem \ref{th1}.
If $\Lambda$ accumulates
$\sigma$ by one side then
$\Lambda$ is transitive.
\end{clly}

\begin{proof}
Let $p,q\in \Lambda$ be periodic.
Suppose that none of the alternatives of Proposition \ref{co0}
hold for $p,q$.
By Theorem \ref{th1'} and Lemma \ref{acc->int} we can assume
$$
W^u_X(p)\cap W^{s,+}_X(\sigma)\neq\emptyset\,\,\,\,
\mbox{and}
\,\,\,\,
W^u_X(q)\cap W^{s,-}_X(\sigma)\neq\emptyset.
$$
Since $\Lambda$ is attracting we have that
$W^u_X(p)\cup W^u_X(q)\subset \Lambda$, and so,
$\Lambda$ intersects both $W^{s,+}_X(\sigma)$ and
$W^{s,-}_X(\sigma)$.
It follows that $\Lambda$ accumulates
$\sigma$ in both sides
contradicting the assumption.
We conclude that one of the two alternatives
in Proposition \ref{co0} holds
$\forall p,q\in Per_X(\Lambda)$.
It follows from this proposition that
$\Lambda$ is transitive and the result follows.
\end{proof}

\begin{prop}
\label{co3}
Let $\Lambda$ be a singular-hyperbolic set
of $X$ satisfying (1)-(4)
of Theorem \ref{th1}.
Suppose that there is $a\in W^u_X(\sigma)\setminus \{\sigma\}$
such that
$\sigma\in \omega_X(a)$.
Then $\Lambda$ is transitive.
\end{prop}

\begin{proof}
Without loss of generality we
can assume that there is
$z\in \omega_X(a)\cap W^+_X(\sigma)$.
If
$W^u_X(q)\cap W^{s,-}_X(\sigma)=\emptyset$
for all $q\in Per_X(\Lambda)$, then
Theorem \ref{th1'} and Lemma
\ref{acc->int} imply
$W^u_X(q)\cap W^{s,+}_X(\sigma)\neq\emptyset$
for all $q\in Per_X(\Lambda)$.
Then $\Lambda$ would be transitive by Proposition
\ref{co0} since the alternative
(1) holds $\forall p,q\in Per_X(\Lambda)$.
So, we can assume that there is
$q\in Per_X(\Lambda)$ such that
$W^u_X(q)\cap W^{s,-}_X(\sigma)\neq \emptyset$.
It follows from the dominating
condition of the singular-hyperbolic splitting of
$\Lambda$ that
the intersection $W^u_X(q)\cap W^{s,-}_X(\sigma)$
is transversal.
This allows us to choose
a point in
$W^u_X(q)$ arbitrarily close to
$W^{s,-}_X(\sigma)$ in the side accumulating
$a$.
Then, as $\Lambda$ is attracting
and satisfies (H1),
it is not difficult to construct
a sequence $p_n\in Per_X(\Lambda)$
converging to
$z\in W^{s,+}_X(\sigma)$,
and $p_n'$ in the orbit of $p_n$
converging to some point in $W^{s,-}_X(\sigma)$.
Suppose by contradiction that $\Lambda$ is not transitive.
Then Proposition \ref{l2}
would imply
$$
W^u_X(p_n')\cap W^{s,+}_X(\sigma)=\emptyset
\,\,\,\,\,\,\mbox{and}
\,\,\,\,\,\,\,W^u_X(p_n)\cap W^{s,-}_X(\sigma)=\emptyset
$$
for $n$ large.
But
$W^u_X(p_n)=W^u_X(p_n')$
since $p_n'$ and $p_n$ are in the same orbit of $X$.
So
$W^u_X(p_n)\cap (W^{s,+}_X(\sigma)\cup W^{s,-}_X(\sigma)
)=\emptyset$.
By Theorem \ref{th1'}
one has
$$
W^u_X(p_n)\cap W^s_X(\sigma)=\emptyset,
$$
a contradiction by Lemma \ref{acc->int}
since $Sing_X(\Lambda)=\{\sigma\}$.
We conclude that $\Lambda$ is transitive
and the proof follows.
\end{proof}

{\flushleft{\bf Proof of Theorem \ref{th1}: }}
Let $\Lambda$ be a singular-hyperbolic set
of a three-dimensional flow $X$
satisfying (1)-(5) of
Theorem \ref{th1}. In particular,
$\Lambda$ has a unique singularity
$\sigma$.
We assume by contradiction that there
is $a\in W^u_X(\sigma)\setminus \{\sigma\}$
such that $\omega_X(a)$ is
{\em not} a periodic orbit.
The contradiction will follow from
Proposition \ref{co0} once we prove that
if $p,q\in Per_X(\Lambda)$ then one of the
alternatives hold:

\begin{enumerate}
\item
$W^u_X(p)\cap W^{s,+}_X(\sigma)\neq\emptyset$
and $W^u_X(q)\cap W^{s,+}_X(\sigma)\neq\emptyset$
\item $W^u_X(p)\cap W^{s,-}_X(\sigma)\neq\emptyset$
and $W^u_X(q)\cap W^{s,-}_X(\sigma)\neq\emptyset$.
\end{enumerate}

Indeed,
by Theorem \ref{th1'} and Lemma \ref{acc->int} we
can assume that
\begin{equation}
\label{estrela}
W^u_X(p)\cap W^{s,+}_X(\sigma)\neq\emptyset\,\,\,\,
\mbox{and}
\,\,\,\,
W^u_X(q)\cap W^{s,-}_X(\sigma)\neq\emptyset.
\end{equation}

\begin{figure}[htv] 
\centerline{
\psfig{figure=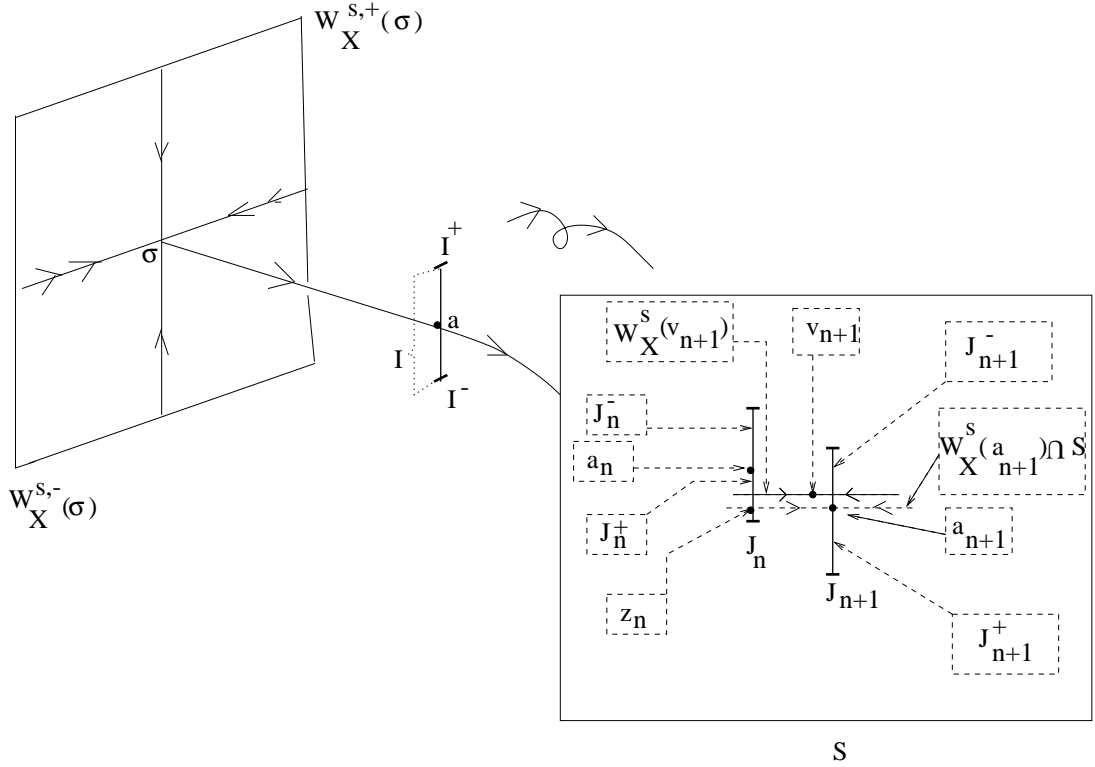,height=4in}}
\caption{\label{f.2} The case $z_n\in J_{n}^+$.}
\end{figure}

By Lemma \ref{lema*}
there is an open interval
$I=I_a$, contained in a suitable
cross-section $\Sigma=\Sigma_a$ of $X$ containing $a$,
such that (1)-(2) of the lemma hold.
See Figure \ref{f.2}.
In particular,
if $I^+, I^-$
are the connected components of $I\setminus \{a\}$ then
$I^+\subset W^u_X(p)$ and
$I^-\subset W^u_X(q)$ (see Figure \ref{f.2}).
Furhermore the tangent direction
of $I$ is contained in $E^c$,
i.e.
$T_xI=E^c_x\cap T_x\Sigma$ for all
$x\in I$.

By Proposition \ref{co3} one has
$\sigma\notin \omega_X(a)$
since
$\Lambda$ is not transitive
(hypothesis (5) of Theorem \ref{th1}).
Then, $H=\omega_X(a)$ is a hyperbolic
saddle-type set by Lemma \ref{lhyp}.
As in \cite{M1} one proves
that $H$ is one-dimensional.
Then we can apply
the Bowen's Theory
of hyperbolic one-dimensional sets \cite{Bo}.
More precisely, there are $\epsilon_0$ small
and a pairwise disjoint family ${\cal S}=\{
S_1,\cdots ,S_r\}$ of local cross sections of time
$\epsilon_0$. Each local cross section
in the family has diameter
less than $\alpha>0$. In addition,
$H=X_{[-\alpha,0]}(H\cap int({\cal S}'))=
X_{[0,\alpha]}(H\cap int({\cal S}'))$,
where
$
{\cal S}'=\cup S_i$
and $int({\cal S}')$ denotes the interior
of ${\cal S}'$.
On the other hand,
$I\subset
\Lambda$ since $\Lambda$
is attracting.
Recall that the tangent direction of $I$
is contained in $E^c$. Moreover,
$E^c$ is volume expanding and
$H$ is non-singular
($Sing_X(\Lambda)=\{\sigma\}$).
Then, the Poincar\'e map
induced by $X$ on ${\cal S}'$
is expanding along $I$.
As
in \cite[p. 371]{MP1}
we can find $\delta>0$ and a open arc sequence
$J_n\subset {\cal S}'$ in the positive orbit
of $I$ satisfying :

\begin{description}
\item{(a)}
the lenght of $J_n$ is
uniformly greather than $\delta$;
\item{(b)}
there is $a_n$ {\em in the positive orbit of $a$}
contained in the interior of
$J_n$.
\end{description}

We can fix $S=S_i\in {\cal S}$
in order to assume that
$J_n\subset S$ for
every $n$.
Let $x\in S$ be a limit point of the sequence
$a_n$.
Then $x\in H\cap int({\cal S}')$.
As the interval $I$ is tangent to the central direction
$E^c$ we have that the interval sequence $J_n$
converges in the $C^1$ topology to an interval
$J\subset W^u_X(x)$.
Observe that $J$ is not trivial by (a).

If $a_n\in W^{s}_X(x)$
for $n$ large we would obtain that $x$
is periodic \cite[Lemma 5.6]{MP1}
which contradicts the assumption
that $\omega_X(a)$ is not periodic.
Then $a_n\notin W^s_X(x)$ for every $n$.
As $J_n\to J$ in the $C^1$ topology
and $\Lambda$ has strong stable manifolds
of uniformly size, we have that there is
$n$ large and
$$
z_n\in
(W^s_X(a_{n+1})\cap S)\cap
(J_{n}\setminus \{a_n\}).
$$

For every $n$
let $J_n^+$ and $J^-_n$ denote the two connected components
of $J_n\setminus \{a_n\}$ in a way that
$J^+_n$ is in the positive orbit of $I^+$ and
$J^-_n$ is in the positive orbit
of $I^-$.
The proof breaks in two cases.

\vspace{5pt}

\noindent
{\em Case I: }$z_n\in J_{n}^+$ (see Figure \ref{f.2}).
By (H1) there is
$v_{n+1}\in Per_X(\Lambda)\cap S$ close to $a_{n+1}$
such that
$$
W^s_X(v_{n+1})\cap J^+_{n}\neq\emptyset
\,\,\,\,
\mbox{and}
\,\,\,\,
W^s_X(v_{n+1})\cap J^-_{n+1}\neq\emptyset.
$$
As $v_{n+1}$ is periodic, close to $a_{n+1}$
and $W^s_X(v_{n+1})\cap J^-_{n+1}\neq\emptyset$
we conclude that
the orbit of $v_{n+1}$ passes
arbitrarily close to some
point in $W^{s,-}_X(\sigma)$.
Then,
$W^u_X(v_{n+1})\cap W^{s,+}_X(\sigma)=\emptyset$
by Proposition \ref{l2} since $\Lambda$ is not transitive.
It follows from Theorem \ref{th1'} that
\begin{equation}
\label{incli}
W^u_X(v_{n+1})\cap W^{s,-}_X(\sigma)
\neq\emptyset.
\end{equation}
But
$W^s_X(v_{n+1})$ intersects
$J^+_{n}$ transversally.
Applying the Inclination Lemma
and (\ref{incli})
we conclude that the positive orbit
of $J^+_{n}$ intersects $W^{s,-}_X(\sigma)$.
As such a positive orbit
is in $W^u_X(p)$ we conclude
that
$$
W^u_X(p)\cap W^{s,-}_X(\sigma)\neq\emptyset.
$$
So, in this case, (\ref{estrela}) implies
$$
W^u_X(p)\cap W^{s,-}_X(\sigma)\neq\emptyset\,\,\,\,
\mbox{and}
\,\,\,\,
W^u_X(q)\cap W^{s,-}_X(\sigma)\neq\emptyset
$$
proving that alternative (1) holds.

\vspace{5pt}

\noindent
{\em Case II: }
$z_n\in J_{n}^-$.
As in the previous case there is
$w_{n+1}\in Per_X(\Lambda)\cap S$ close to $a_{n+1}$
such that
$$
W^s_X(w_{n+1})\cap J^-_{n}\neq\emptyset
\,\,\,\,
\mbox{and}
\,\,\,\,
W^s_X(w_{n+1})\cap J^+_{n+1}\neq\emptyset.
$$
As before the orbit of $w_{n+1}$ passes close
to $W^{s,+}_X(\sigma)$.
Since $\Lambda$ is not transitive
we have by Proposition \ref{l2} and Theorem \ref{th1'}
as before that
$$
W^u_X(w_{n+1})\cap W^{s,+}_X(\sigma)
\neq\emptyset.
$$
Again
$W^s_X(w_{n+1})$ intersects
$J^-_{n}$ transversally.
Then, the Inclination Lemma as before
implies that the positive orbit
of $J^-_{n}$ intersects $W^{s,+}_X(\sigma)$.
As such a positive orbit
is in $W^u_X(q)$ we conclude
that
$$
W^u_X(q)\cap W^{s,+}_X(\sigma)\neq\emptyset.
$$
So, in this case, applying (\ref{estrela})
one has
$$
W^u_X(p)\cap W^{s,+}_X(\sigma)\neq\emptyset\,\,\,\,
\mbox{and}
\,\,\,\,
W^u_X(q)\cap W^{s,+}_X(\sigma)\neq\emptyset
$$
proving that alternative (2) holds.
\qed

{\flushleft{\bf Proof of Theorem \ref{generico}: }}
By the Kupka-Smale Theorem \cite{dMP}
the set of Kupka-Smale
flows ${\cal R}$ is
residual in ${\cal X}^r$, $\forall r\geq 1$.
If $X$ is Kupka-Smale
then the unstable manifold
of any Lorenz-like
singularity $\sigma$ of $X$ cannot intersect the stable
manifold of any hyperbolic saddle periodic orbit of
$X$
(because such a intersection is not transversal).
So, there is no
$a\in W^u_X(\sigma)\setminus \{\sigma\}$
such that $\omega_X(a)$ is a periodic orbit.
Then the result follows from Theorem \ref{th1}.
\qed




\section{Proof of Theorem \ref{thA'} and Corollary \ref{galinha}}

Theorem \ref{thA'} will be a direct consequence
of Corollary \ref{cara} and
the following result.

\begin{thm}
\label{homoclinic class}
Let $\Lambda$ be a singular-hyperbolic set
of a $C^r$ flow $X$ on a closed three-manifold,
$r\geq 1$.
Suppose that the properties below hold.

\begin{enumerate}
\item
$\Lambda$ is connected.
\item
$\Lambda$ is attracting.
\item
The closed orbits contained in $\Lambda$ are dense in $\Lambda$.
\item
$\Lambda$ has a unique singularity $\sigma$.
\item
$\Lambda$ is not transitive.
\end{enumerate}
Then,
$H^+_X(\sigma,\Lambda)$
and $H^-_X(\sigma,\Lambda)$
are homoclinic classes of $X$.
\end{thm}

Before the proof we state some preliminary results.
First of all
we shall assume, along this section, that
$\Lambda$ and $X$ satisfies
(1)-(5) of Theorem \ref{homoclinic class}.
It follows from Theorem \ref{th1} that
$\omega_X(a)$ is a periodic orbit for
every $a\in W^u_X(\sigma)\setminus \{\sigma\}$.
Hereafter we fix
$a\in W^u_X(\sigma)\setminus \{\sigma\}$,
and $O_0$ periodic with $\omega_X(a)=O_0$.
We also fix $p_0\in O_0$.
It is clear that $O_0$ is saddle-type.
Then, as $X$ is three-dimensional,
the Poincar\'e map
induced by $X$ in a
suitable cross-section at $p_0$ has only one
eigenvalue with modulus $>1$.
We shall call such a eigenvalue as
the {\em expanding eigenvalue} of
$p_0$.
The lemma below will be used
only in Definition \ref{u,+-}.

\begin{figure}[htv] 
\centerline{
\psfig{figure=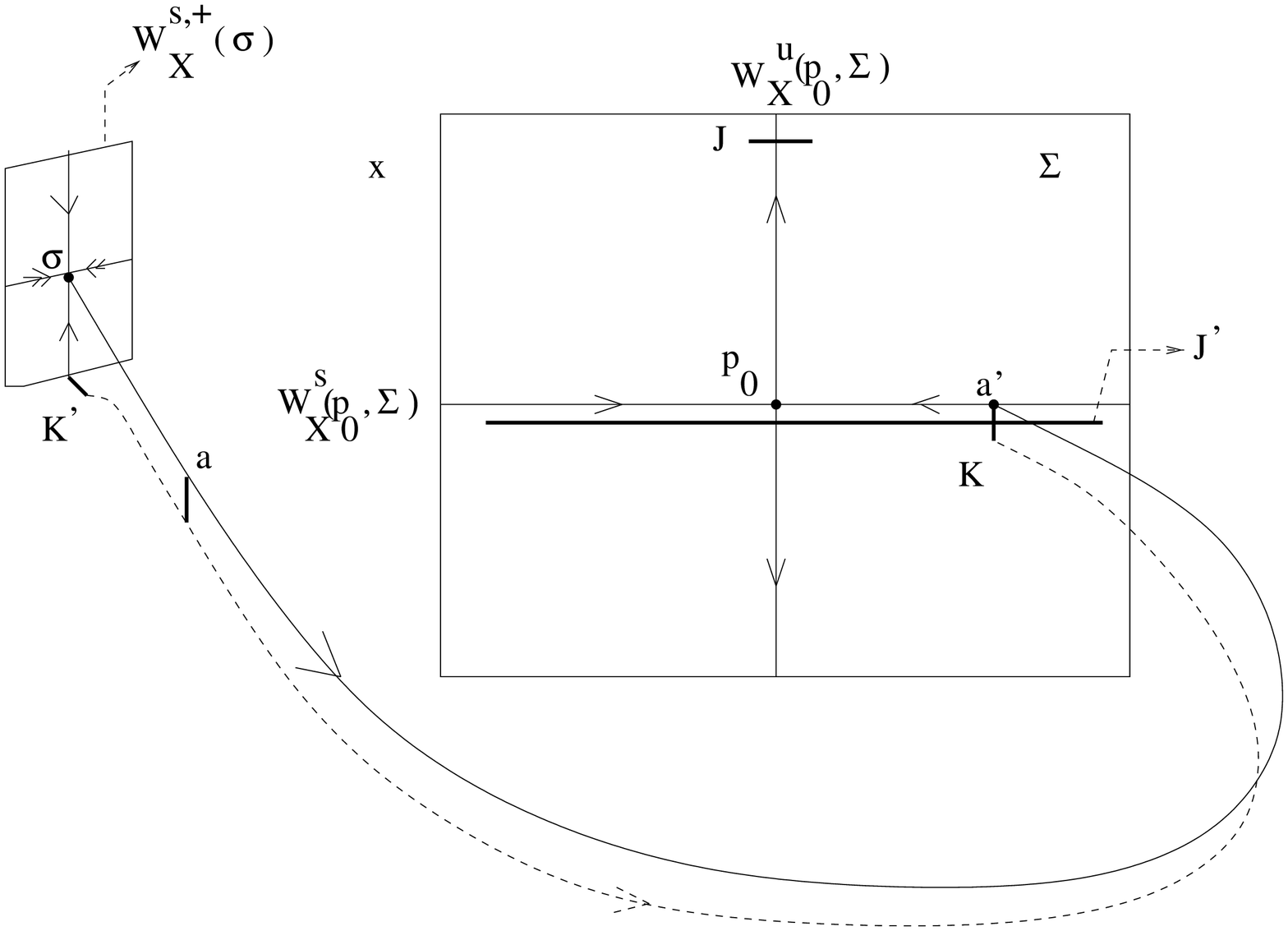,height=3in}}
\caption{\label{f.4} Positive expanding eigenvalues.}
\end{figure}  

\begin{lemma}
\label{c.autovalorreais}
The expanding eigenvalue
of $p_0$ is positive.
\end{lemma}

\begin{proof}
Suppose by contradiction that the
expanding eigenvalue of
$p_0$ is negative.
Fix a cross-section $\Sigma$ of $X$
containing $p_0$.
This section defines a Poincar\'e map
$\Pi:Dom(\Pi)\subset \Sigma\to \Sigma$.
By definition the eigenvalues
of the derivative $D\Pi(p_0)$
are precisely the ones of $p_0$.
It follows that $D\Pi(p_0)$
has a negative eigenvalue.
As $p_0\in Per_X(\Lambda)$,
Theorem \ref{th1'} and
Lemma \ref{acc->int}
imply
that $W^u_X(p_0)$ intersects
$W^{s,+}_X(\sigma)$ or
$W^{s,-}_X(\sigma)$.
We shall assume the first case
since the proof for the second one is similar.
Under such conditions we shall prove that
$W^u_X(p)\cap W^{s,+}_X(\sigma)\neq\emptyset$ for all
$p\in Per_X(\Lambda)$. For this we proceed as follows.
Let $p\in Per_X(\Lambda)$ be fixed.
Again by Theorem \ref{th1'} and Lemma \ref{acc->int}
we have that
$W^u_X(p)$ intersects
$W^{s,+}_X(\sigma)$ or
$W^{s,-}_X(\sigma)$.
In the first case we are done. So, we can assume that
$W^u_X(p)\cap W^{s,-}_X(\sigma)\neq \emptyset$.
It follows from
Lemma \ref{sigma} that
$W^u_X(p)$ contains an interval
$K'=J_p$ as indicated in Figure
\ref{f.4}.
By the definition of $p_0$ there is
$a\in W^u_X(\sigma)\setminus \{\sigma\}$
such that $\omega_X(a)$ is the full orbit of $p_0$
(in particular $a\in W^s_X(p)$).
We can assume that $K'$ points in the side of
$a$ as indicated in the figure.
Then,
as $a\in W^s_X(p_0)$, the positive iterate
of $K'$ first passes through $a$ and then goes to
$\Sigma$ yielding the
interval $K$ in Figure \ref{f.4}.
Note that the closure
of $K$ intersects $W^s_X(p_0,\Sigma)$
transversally at $a'$ as in Figure \ref{f.4}
(observe that $a'$ is in the positive orbit of $a$).
On the other hand, since $W^u_X(p_0)\cap W^{s,+}_X(\sigma)
\neq\emptyset$, we have that there is
an interval
$J\subset W^{s,+}_X(\sigma)\cap\Sigma$ intersecting
$W^u_X(p_0,\Sigma)$ transversally as indicated in Figure
\ref{f.4}.
Since the expanding eigenvalue
of $D\Pi(p_0)$
is negative we have by the Inclination Lemma that
the negative $\Pi$-iterates of $J$
accumulates $W^s_X(p_0,\Sigma)$
{\em in both sides} of
$\Sigma\setminus W^s_X(p_0,\Sigma)$.
In particular,
as $K$ intersects
$W^s_X(p_0,\Sigma)$ transversally at $a'$,
there is a negative iterate
$J'$ of $J$ intersecting
$K$ as indicated in Figure \ref{f.4}.
It follows from the invariance of the
stable and unstable manifolds
that such a intersection
belongs to $W^u_X(p)\cap W^{s,+}_X(\sigma)$
since $J\subset W^{s,+}_X(\sigma)$ and
$K\subset W^u_X(p)$.
This proves
$W^u_X(p)\cap W^{s,+}_X(\sigma)\neq\emptyset$ for all
$p\in Per_X(\Lambda)$.
So, the alternative
(1) in Proposition
\ref{co0} holds
for all $p,q\in Per_X(\Lambda)$
proving that
$\Lambda$ is transitive, a contradiction.
This contradiction shows that
the expanding eigenvalue
of $p_0$ is positive.
\end{proof}

\begin{figure}[htv] 
\centerline{
\psfig{figure=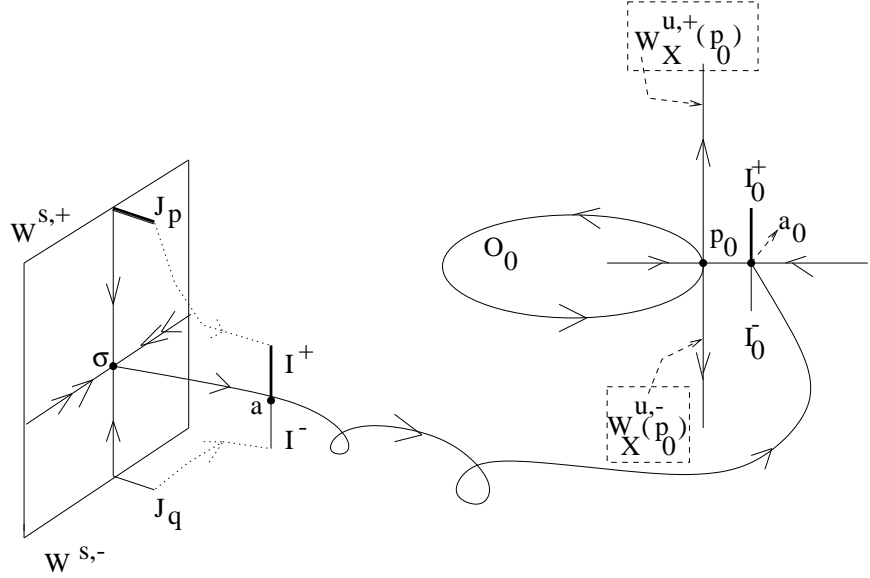,height=3in}}
\caption{\label{f.5} Definition of $W^{u,+}_X(p_0)$ and
$W^{u,-}_X(p_0)$.}
\end{figure}

Let us explain how
Lemma \ref{c.autovalorreais} applies
in the proof of Theorem \ref{homoclinic class}.
By
Proposition \ref{co0} and
(5) of Theorem \ref{homoclinic class}
there are $p,q\in Per_X(\Lambda)$ such that
$$
W^u_X(p)\cap W^{s,+}_X(\sigma)\neq\emptyset\,\,\,\,
\mbox{and}
\,\,\,\,
W^u_X(q)\cap W^{s,-}_X(\sigma)\neq\emptyset.
$$
By Lemma \ref{lema*}
there is
an interval
$I=I_a$, contained in a suitable cross-section of $X$ and containing $a$, such that
if $I^+,I^-$ are the connected components
of $I\setminus \{a\}$ then
$I^+\subset W^u_X(p$ and $I^-\subset W^u_X(q)$.
In addition $I$ is tangent to the central direction $E^c$ of $\Lambda$ (see Figure \ref{f.5}).
As $a\in W^s_X(p_0)$ and $I$ is tangent to
$E^c$ we have that the flow
of $X$ carries $I$ to an interval $I'$
transverse to $W^s_X(p_0)$ at $a_0$
as indicated in Figure
\ref{f.5}.
Note that the flow carries $I^+$ and
$I^-$ into $I^+_0$ and $I^-_0$ respectively.
On the other hand, Lemma \ref{c.autovalorreais}
implies that
$W^u_X(p_0)$ is a cylinder with generating curve
$O_0$.
In particular, $O_0$ separates $W^u_X(p_0)$
in two connected components.
We denote such components by
$W^{u,+}_X(p_0)$ and
$W^{u,-}_X(p_0)$
according to the following definition.

\begin{defi}
\label{u,+-}
We denote by
$W^{u,+}_X(p_0)$ the connected component
of $W^{u}_X(p_0)\setminus O_0$
which is accumulated
(via Inclination Lemma) by the positive orbit
of $I^+_0$.
We denote by
$W^{u,-}_X(p_0)$ the connected component
of $W^{u}_X(p_0)\setminus O_0$
which is accumulated by the positive orbit
of $I^-_0$.
\end{defi}

\begin{lemma}
\label{independence}
Definition \ref{u,+-}
does not depend on $p,q,J_p,J_q$.
\end{lemma}

\begin{proof}
This lemma is a
straightforward consequence
of the Inclination
Lemma, the Strong $\lambda$-lemma
and Lemma \ref{c.autovalorreais}.
The details are left to the reader.
\end{proof}

The next propositions give some properties
of $W^{u,+}_X(p_0)$ and $W^{u,-}_X(p_0)$.
In particular, Theorem \ref{homoclinic class}
will follow from
Propositions \ref{homoclinic3}
and \ref{cl=homoclinic}.

\begin{prop}
\label{preliminar}
$W^{u,+}_X(p_0)\cap W^{s,+}_X(\sigma)\neq\emptyset$
and $W^{u,+}_X(p_0)\cap W^{s,-}_X(\sigma)=\emptyset$.
Similarly interchanging the roles of $+$ and $-$.
\end{prop}

\begin{proof}
For simplicity denote
$W=W^{u,+}_X(p_0)$.
First we prove
$W\cap W^{s,-}_X(\sigma)=\emptyset$.
Suppose by contradiction that
$W\cap W^{s,-}_X(\sigma)\neq\emptyset$.
As this last intersection is transversal
there is an interval $J\subset
W^{s,-}_X(\sigma)$ intersecting
$W$ transversally.
Now, fix a cross-section $\Sigma=\Sigma^+$ as in
Lemma \ref{sigma} and let
$p\in Per_X(\Lambda)$ be such that
$W^u_X(p)\cap W^{s,+}_X(\sigma)\neq\emptyset$.
Then, there
is an small interval $I\subset
W^u_X(p)\cap \Sigma$ transversal to
$\Sigma\cap W^{s,+}_X(\sigma)$.
By the definition
of $W=W^{u,+}_X(p_0)$
(Definition \ref{u,+-})
we have that the positive orbit
of $I$ accumulates
on $W$.
As $J$ is transversal to
$W$
the Inclination Lemma implies that
the positive orbit of $I$ intersects
$J$. This proves
$W^u_X(p)\cap W^{s,-}_X(\sigma)\neq\emptyset$
for all $p\in Per_X(\Lambda)$.
It would follow that alternative
(2) of Proposition \ref{co0}
holds $\forall p,q$ contradicting
(5) of Theorem \ref{homoclinic class}.
This contradiction proves
$W\cap W^{s,-}_X(\sigma)=\emptyset$
as desired.
Now suppose by contradiction that
$W\cap W^{s,+}_X(\sigma)=\emptyset$.
As $W\cap W^{s,-}_X(\sigma)=\emptyset$
we would obtain from Theorem \ref{th1'}
that $W\cap W^s_X(\sigma)=\emptyset$.
But the denseness of the periodic orbits
together with the Inclination Lemma
imply $W\cap W^s_X(\sigma)\neq\emptyset$.
This contradiction proves
$W\cap W^{s,+}_X(\sigma)\neq\emptyset$
and the result follows.
\end{proof}

\begin{prop}
\label{homoclinic3}
$H^+_X(\sigma,\Lambda)=Cl(W^{u,+}_X(p_0))$.
Similarly replacing $+$ by $-$.
\end{prop}

\begin{proof}
Fix $q\in P^+_X(\sigma,\Lambda)$,
i.e.
$W^u_X(q)\cap W^{s,+}_X(\sigma)\neq\emptyset$.
Note that $W^{u,+}_X(p_0)\cap
W^{s,+}_X(\sigma)\neq\emptyset$ by Lemma
\ref{preliminar}.
Then, as in the proof of
Lemma \ref{l1}, one can prove that
$W^{u,+}_X(p_0)$ accumulates on
$q$
(in this case
$p=p_0$ and $z_q\in W^{u,+}_X(p_0)$
at Figure \ref{f.4}).
This proves
$
H^+_X(\sigma,\Lambda)\subset
Cl(W^{u,+}_X(p_0))
$.
Conversely
let $x\in W^{u,+}_X(p_0)$ be fixed.
Since the periodic orbits are dense
in $\Lambda$ and
$W^{u,+}_X(p_0)\subset \Lambda$
there is
$z\in Per_X(\Lambda)$ nearby $x$.
Choosing $z$ close to $x$
we assure
$W^s_X(z)\cap W^{u,+}_X(p_0)\neq\emptyset$
because stable manifolds have size
uniformly bounded away from zero.
If $W^u_X(z)\cap W^{s,-}_X(\sigma)\neq \emptyset$ then
the Inclination Lemma
and
$W^s_X(z)\cap W^{u,+}_X(p_0)\neq\emptyset$
would imply
$W^u_X(p_0)\cap W^{s,-}_X(\sigma)\neq
\emptyset$.
This contradicts Proposition \ref{preliminar}
and so
$W^u_X(z)\cap W^{s,-}_X(\sigma)= \emptyset$.
By Theorem \ref{th1'} and Lemma
\ref{acc->int} we conclude that
$z\in P^+_X(\sigma,\Lambda)$
and then
$x\in H^+_X(\sigma,\Lambda)$
by Definition \ref{H+}.
The lemma is proved.
\end{proof}

\begin{prop}
\label{dense II}
If $z\in Per_X(\Lambda)$ and
$W^s_X(z)\cap W^{u,+}_X(p_0)\neq\emptyset$,
then
$Cl(W^s_X(z)\cap W^{u,+}_X(p_0))=
Cl(W^{u,+}_X(p_0))$.
Similarly replacing $+$ by $-$.
\end{prop}
\begin{proof}
First we show
$Cl(W^{u,+}_X(p_0)\cap W^{s,+}_X(\sigma))=Cl(W^{u,+}_X(p_0))$.
By Proposition \ref{homoclinic3} one has
$Cl(W^{u,+}_X(p_0))=H^+_X(\sigma,\Lambda)$.
Fix $x\in W^{u,+}_X(p_0)$.
Since the periodic orbits are dense
in $\Lambda$ there is
$w\in Per_X(\Lambda)$ close to $x$.
Suppose that $W^u_X(w)\cap W^{s,+}_X(\sigma)
=\emptyset$.
Then $W^u_X(w)\cap W^{s,-}_X(\sigma)
=\emptyset$ by Theorem \ref{th1'} and Lemma
\ref{acc->int}.
It follows from the Inclination Lemma
and Proposition \ref{preliminar} that
$W^{u,+}_X(p_0)\cap W^{s,-}_X(\sigma)\neq\emptyset$
contradicting Proposition \ref{preliminar}.
This proves 
$W^u_X(w)\cap W^{s,+}_X(\sigma)
\neq\emptyset$.
Note that $W^s_X(w)\cap W^{u,+}_X(p_0)\neq \emptyset$
is close to $x$ (as $w\to x$).
As
$W^u_X(w)\cap W^{s,+}_X(\sigma)
\neq\emptyset$ is transversal we
can apply the Inclination Lemma
in order to find
a transverse intersection
$W^{u,+}_X(p_0)\cap W^{s,+}_X(\sigma)$ close to $x$. This proves
$Cl(W^{u,+}_X(p_0)\cap W^{s,+}_X(\sigma))=Cl(W^{u,+}_X(p_0))$.
Now we prove
$Cl(W^s_X(z)\cap W^{u,+}_X(p_0))=
Cl(W^{u,+}_X(p_0))$.
Choose
$x\in W^{u,+}_X(p_0)$.
Then there is an interval
$I_x\subset W^{u,+}_X(p_0)$ arbitrarily close
to $x$ such that
$I_x\cap W^{s,+}_X(\sigma)\neq\emptyset$.
The positive orbit
of $I_x$ first passes through
$a$ and, afterward, it accumulates
on $W^{u,+}_X(p_0)$
(see Figure \ref{f.5}).
But $W^s_X(z)\cap W^{u,+}_X(p_0)\neq\emptyset$ by assumption. As this intersection is transversal
the Inclination Lemma implies
that the positive orbit of
$I_x$ intersects $W^s_X(z)$.
By taking the backward flow of the last intersection
we get $W^s_X(z)\cap I_x\neq\emptyset$.
This proves the lemma.
\end{proof}

\vspace{10pt}

\begin{figure}[htv] 
\centerline{
\psfig{figure=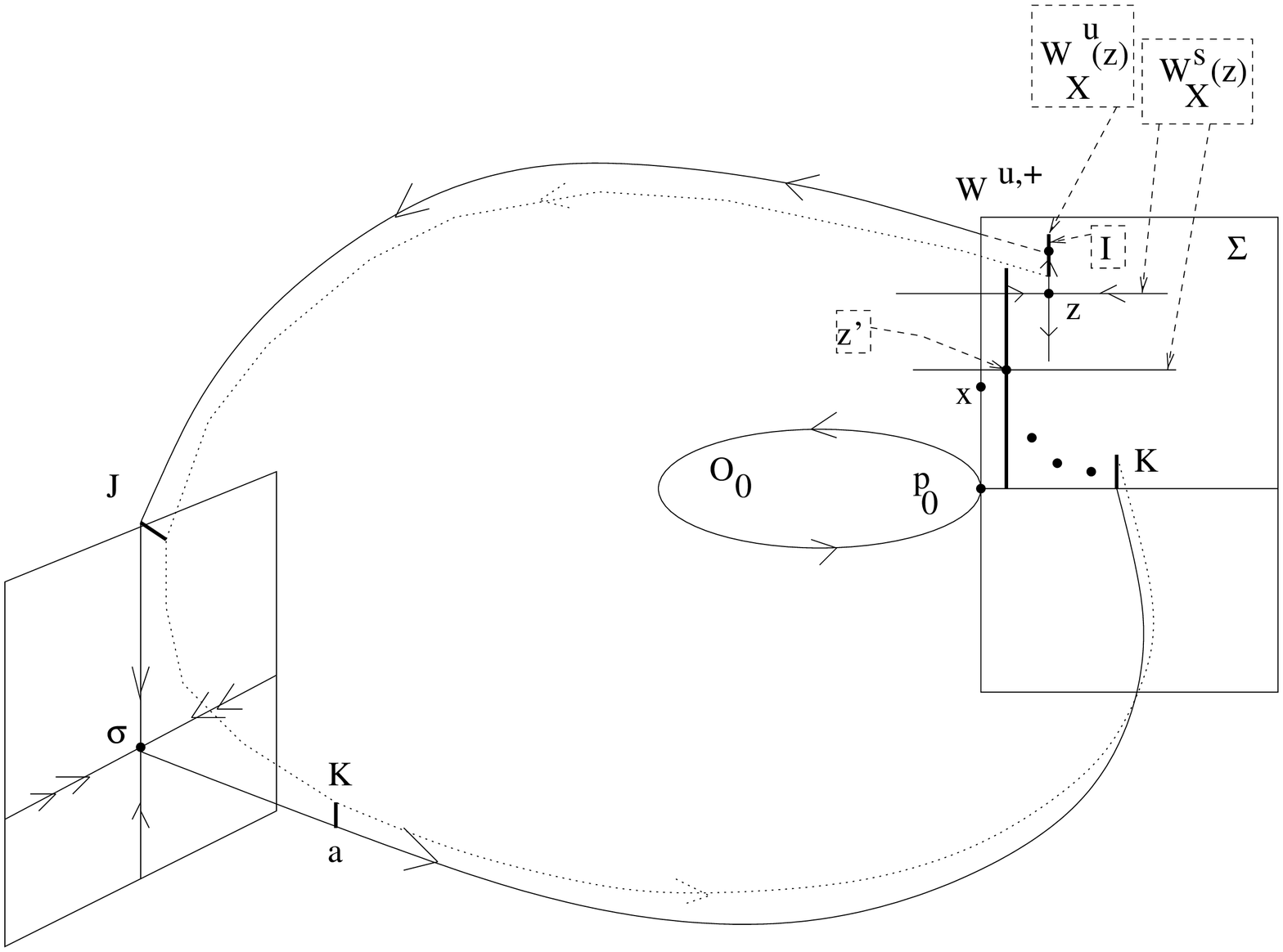,height=3in}}
\caption{\label{f.6} $Cl(W^{u,+}_X(p_0))$
is a homoclinic class.}
\end{figure}

Given $z\in Per_X(\Lambda)$ we denote by
$H_X(z)$ the homoclinic class associated to $z$.

\begin{prop}
\label{cl=homoclinic}
If $z\in Per_X(\Lambda)$ is close to a point in
$W^{u,+}_X(p_0)$, then
$H_X(z)=Cl(W^{u,+}_X(p_0))$.
Similarly replacing $+$ by $-$.
\end{prop}

\begin{proof}
Let $z\in Per_X(\Lambda)$ be a point
close to one in $W^{u,+}_X(p_0)$.
It follows from the continuity
of the stable manifolds that
$
W^s_X(z)\cap W^{u,+}_X(p_0)\neq\emptyset$.
We claim that
$H_X(z)=Cl(W^{u,+}_X(p_0))$.
Indeed, by Proposition \ref{preliminar}
one has $W^{u,+}_X(p_0)\cap W^{s,-}_X(\sigma)=
\emptyset$. This equality and the Inclination Lemma
imply $W^u_X(z)\cap W^s_X(\sigma)\subset
W^{s,+}_X(\sigma)$.
By
Proposition \ref{dense II} we have that
$W^s_X(z)\cap W^{u,+}_X(p_0)$ is dense
in $W^{u,+}_X(p_0)$ since
$W^s_X(z)\cap W^{u,+}_X(p_0)\neq\emptyset$.
Let $\Sigma$ be a cross-section containing $p_0$
as in Figure \ref{f.6} and
fix $x\in W^{u,+}_X(p_0)$.
We can assume that $x,z\in \Sigma$
are as in Figure \ref{f.6}.
Since $W^u_X(z)\cap W^s_X(\sigma)\neq\emptyset$
and
$W^u_X(z)\cap W^s_X(\sigma)\subset W^{s,+}_X(\sigma)$
there is an interval
$I\subset W^u_X(z)$ intersecting
$W^{s,+}_X(\sigma)$.
Then, the positive orbit of $I$ yields the interval
$J$ in the figure.
In addition,
the positive orbit of $J$ yields the interval
$K$ in that figure.
Note that the positive orbit of $K$
{\em accumulates} $W^{u,+}_X(p_0)$
as in Figure \ref{f.6}
(recall Definition
\ref{u,+-}).
As $W^s_X(z)\cap W^{u,+}_X(p_0)$
is dense in $W^{u,+}_X(p_0)$
and $x\in W^{u,+}_X(p_0)$
we have that
$W^s_X(z)$ passes close to $x$ as indicated in the figure.
The Inclination Lemma applied to the positive orbit
of $K$ yields a homoclinic point
$z'$ associated to $z$ which is close to $x$.
This proves that $x\in H_X(z)$, and so, $Cl(W^{u,+}_X(p_0))\subset H_X(z)$.
The reversed inclusion
is a direct consequence of the Inclination Lemma
applied to $W^s_X(z)\cap W^{u,+}_X(p_0)\neq\emptyset$.
We conclude that $Cl(W^{u,+}_X(p_0))=
H_X(z)$ proving the result.
\end{proof}

{\flushleft{\bf Proof of Theorem
\ref{homoclinic class}: }}
Let $\Lambda$ be a singular-hyperbolic set
of $X$ satisfying (1)-(5) of
Theorem \ref{homoclinic class}.
Then we can apply the results in this section.
To prove that $H^+_X(\sigma,\Lambda)$ is
a homoclinic classs it suffices
by Proposition \ref{homoclinic3} to
prove that $Cl(W^{u,+}_X(p_0))$ is a homoclinic
class.
By (3) of Theorem \ref{homoclinic class}
we can choose $z\in Per_X(\Lambda)$
arbitrarily close to
a point in $W^{u,+}_X(p_0)$.
Then
$Cl(W^{u,+}_X(p_0))=H_X(z)$
by Proposition \ref{cl=homoclinic}
and the result follows.
\qed

{\flushleft{\bf Proof of Corollary
\ref{galinha}: }}
Let $\Lambda$ be
a singular-hyperbolic set
of a three-dimensional flow $X$ which is
connected, attracting, have dense closed orbits and only
one  singularity (denoted by $\sigma$).
By Theorem \ref{th1'} we have that
$\sigma$ is Lorenz-like, and so,
$W^u_X(\sigma)\setminus \{\sigma\}$
consists of two regular orbits.
Fix $a,a'$ in each of such orbits.
We shall assume that
$\Lambda$ is not transitive.
Then, $\Lambda$ satisfies the hypotheses
(1)-(5) of Theorem \ref{th1}.
It follows from this theorem that
$\omega_X(a)=O_0$ and
$\omega_X(a')=O_0'$, where
$O_0$ and $O_0'$ are periodic orbits of $X$.
On the other hand, $\Lambda$ satisfies the
hypotheses (1)-(5) of Theorem
\ref{homoclinic class} too.
Then the results of this 
section can be applied.
As before we can fix $p_0'\in O_0'$
and prove
the corresponding version of
the results of this section for $p_0'$.
In particular, we can define
$W^{u,+}_X(p_0'),W^{u,-}_X(p_0')$
as in Definition \ref{u,+-} and prove as
in Proposition
\ref{cl=homoclinic} that if $z\in Per_X(\Lambda)$ is
close to a point in
$W^{u,+}_X(p_0')$ then
$Cl(W^{u,+}_X(p_0'))=
H_X(z)$ (similarly for
$-$).

Now assume by contradiction that $\Lambda$
is not tame, i.e.
it has infinitely many homoclinic classes
$H_X(z_n)$.
We have two possibilities, namely
$\sigma\in Cl(\cup_nH_X(z_n))$ or
$\sigma\notin Cl(\cup_nH_X(z_n))$.
If $\sigma\notin Cl(\cup_nH_X(z_n))$
then $Cl(\cup_nH_X(z_n))$
is a hyperbolic set by Lemma \ref{lhyp}.
This is a contradiction because
hyperbolic sets cannot contain infinitely many
homoclinic classes \cite{N}.
This contradiction proves
$\sigma\in Cl(\cup_nH_X(z_n))$.
Let $\Sigma^+,\Sigma^-$ be the cross-sections of $X$
as in Lemma \ref{sigma} for $\sigma$, i.e.
every orbit of $\Lambda$ passing close to
some point in $W^{s,+}_X(\sigma)$
(resp. $W^{s,-}_X(\sigma)$)
intersects $\Sigma^+$ (resp. $\Sigma^-$).

Because $\sigma\in Cl(\cup_nH_X(z_n))$
the above property of $\Sigma^\pm$ and Theorem \ref{th1'}
imply that
the orbit of $z_n$ can be assumed to intersect
$\Sigma^+\cup \Sigma^-$.
Replacing $z_n$ by one point in its orbit if necessary
we can further assume that
$z_n\in \Sigma^+\cup \Sigma^-$.
Recall that $\Sigma^+,\Sigma^-$ are respectively
equiped with the vertical curves $l^+,l^-$
at Figure \ref{f.-2}.
Note that $(\Sigma^+\setminus l^+)\cup
(\Sigma^-\setminus l^-)$ is formed by
four connected components.
The forward orbits starting in
two of these components passes close to
$a$ (as indicated in Figure \ref{f.-2})
while the forward orbits starting in the
remaining ones passes close to $a'$.
Since the orbit of $z_n$ are periodic ($\forall n$)
we have that $z_n\notin l^+\cup l^-$.
Then, as the sequence $z_n$ is infinite
(and $(\Sigma^+\setminus l^+)\cup
(\Sigma^-\setminus l^-)$ has
four connected components only), we can assume that
$z_n$ is in only one of these components.
It follows that the positive orbit of
$z_n$ passes close to either $a$ or $a'$
($\forall n$).
We are going to assume that the orbit
of $z_n$ passes close to $a$ ($\forall n$)
becase the argument for $a'$
is enterely similar
(for this we use the version of the results
of this section for $a'$ and $p_0'$).
As we can observe in Figure \ref{f.5}
the positive orbit
of $a$ intersects the stable manifold close to $p_0$
in the point $a_0'$. The orbit positive of
$z_n$ (which passes close to $a$) will pass
close to $a_0'$ too.
Clearly such a positive orbit
is not in the stable manifold
of $p_0$.
By usual fundamental domain arguments
\cite{dMP} we have that there is
a sequence $z_n'$ in the positive orbit
of $z_n$ converging to
a point in either $W^{u,+}_X(p_0)$ or
$W^{u,-}_X(p_0)$.
In each case we conclude by
Proposition \ref{cl=homoclinic} that
$H_X(z_n')$ is either $Cl(W^{u,+}_X(p_0))$ or
$Cl(W^{u,+}_X(p_0))$ for all $n$ (say).
Because the sequence $z_n'$ is infinite we can
assume that
$H_X(z_n')=Cl(W^{u,+}_X(p_0))$
for all $n$.
Because $z_n'$ is in the orbit of $z_n$ we conclude
that $H_X(z_n)=Cl(W^{u,+}_X(p_0))$
for all $n$.
This would imply that $H_X(z_n)=H_X(z_m)=
Cl(W^{u,+}_X(p_0))$ for all $n,m$ a contradiction.
This contradiction proves the result.
\qed

\vspace{0.2cm}
\noindent C. A. Morales, M. J. Pacifico\\
Instituto de Matem\'atica \\
Universidade Federal do Rio de Janeiro \\
C. P. 68.530, CEP 21.945-970 \\
Rio de Janeiro, R. J. , Brazil \\
e-mail: \verb morales@impa.br, \verb
pacifico@im.ufrj.br

\end{document}